\documentclass[11pt]{article}
\usepackage{amsmath}
\usepackage{amsfonts}
\usepackage{amssymb}
\usepackage{dsfont}
\usepackage{hyperref}
\usepackage{indentfirst}
\usepackage{mathrsfs}

\usepackage[top=1in,bottom=1in,left=1.25in,right=1.25in]{geometry}
\date{}
\allowdisplaybreaks
\begin{document}

\renewcommand{\baselinestretch}{1.2}
\renewcommand{\arraystretch}{1.0}

\title{\bf Multiplier Hopf coquasigroups with faithful integrals}
\author
{
 \textbf{Tao Yang} \footnote{Corresponding author.
 College of Science, Nanjing Agricultural University, Nanjing 210095, China. E-mail: tao.yang@njau.edu.cn}
}
\maketitle

\begin{center}
\begin{minipage}{12.cm}

 \textbf{Abstract}
 Let $A$ be a multiplier Hopf coquasigroup. If the faithful integrals exist, then they are unique up to scalar.
 Furthermore, if  $A$ is of discrete type, then its integral duality $\widehat{A}$ is a Hopf quasigroup, 
 and the biduality $\widehat{\widehat{A}}$ is isomorphic to the original $A$ as multiplier Hopf coquasigroups.
 This biduality theorem also holds for a class of Hopf quasigroups with faithful integrals.
\\

 {\bf Key words} multiplier Hopf coquasigroup, integral, Hopf quasigroup
\\

 {\bf Mathematics Subject Classification}   16T05, 16T99

\end{minipage}
\end{center}
\normalsize

\section{Introduction}
\def\theequation{\thesection.\arabic{equation}}
\setcounter{equation}{0}

 Multiplier Hopf coquasigroups introduced in \cite{Y} gave an answer to the duality of some infinite dimensional Hopf quasigroups.
 This notion generalized the Hopf coquasigroup to a nonunital case.
 
 Here we continue our work on multiplier Hopf coquasigroups. 
 First of all, we consider the uniqueness of the faithful integrals, and get that the faithful integral is unique up to scalar.
 Based on the unique integral, we construct the integral duality of a discrete multiplier Hopf coquasigroup and show that the duality is a Hopf quasigroup introduced in \cite{K, KM}.
 According to the previous main result in \cite{Y}, we find that the biduality $\widehat{\widehat{A}}$ of a discrete multiplier Hopf coquasigroup $A$ is isomorphic to $A$ itself,
 and the biduality $\widehat{\widehat{H}}$ of a Hopf quasigroup $H$ is isomorphic to $H$, 
 i.e., biduality theorem holds for some multiplier Hopf coquasigroups and infinite-dimensional Hopf quasigroup.

 This paper is organized as follows. In section 2, we introduce some notions: multiplier Hopf coquasigroups and integrals, which will be used in the following sections.

 In section 3, we consider the local units and integrals on a multiplier Hopf coquasigroup. 
 We show that the faithful integrals are unique up to scalar, and there exists a modular automorphism.

 In section 4, for a discrete multiplier Hopf coquasigroup $A$ with a faithful integral $\varphi$, we construct the dual $\widehat{A} = \varphi(\cdot A)$, 
 and give the multiplication, comultiplication, counit and antipode to make $\widehat{A}$ a Hopf quasigroup.

 In section 5, we show that the biduality theorem holds for some Hopf quasigroups and discrete multiplier Hopf coquasigroups.
 That is, the bidualities are isomorphic to the original ones.

\section{Preliminaries}
\def\theequation{\thesection.\arabic{equation}}
\setcounter{equation}{0}

 Throughout this paper, all spaces we considered are over a fixed field $k$ (e.g., the complex number field $\mathds{C}$).
 Let $A$ be an (associative) algebra. We do not assume that $A$ has a unit, but we do require that
 the product, seen as a bilinear form, is non-degenerated  (i.e., whenever $a\in A$ and $ab=0$ for all $b\in A$
 or $ba=0$ for all $b\in A$, we must have that $a=0$).
 Then we can consider the multiplier algebra $M(A)$ of $A$.

 Recall from \cite{V94, V98} that $M(A)$ is characterized as the largest algebra with identity containing $A$ as an essential two-sided ideal.
 In particularly, we still have that, whenever $a\in M(A)$ and $ab=0$ for all $b\in A$ or $ba=0$ for all $b\in A$, again $a=0$.
 
 Furthermore, we consider the tensor algebra $A\otimes A$. It is still non-degenerated and we have its multiplier algebra $M(A\otimes A)$.
 There are natural imbeddings
 $$A\otimes A \subseteq M(A)\otimes M(A) \subseteq M(A\otimes A).$$

 In generally, when $A$ has no identity, these two inclusions are stict.
 If $A$ already has an identity, the product is obviously non-degenerate
 and $M(A)=A$ and $M(A\otimes A) = A\otimes A$. 

 Let $A$ and $B$ be non-degenerate algebras, if homomorphism $f: A\longrightarrow M(B)$ is non-degenerated
 (i.e., $f(A)B=B$ and $Bf(A)=B$),
 then has a unique extension to a homomorphism $M(A)\longrightarrow M(B)$, we also denote it $f$.

Recall from \cite{Y} a \emph{multipler Hopf coquasigroup} is a nondegenerate associative algebra $A$ equipped with algebra homomorphisms $\Delta: A\longrightarrow M(A\otimes A)$(coproduct),
 $\varepsilon: A\longrightarrow k$(counit) and a linear map $S: A \longrightarrow A$ (antipode) such that
 \begin{enumerate}
 \item[(1)] $T_{1}(a\otimes b)=\Delta(a)(1\otimes b)$ and $T_{2}(a\otimes b)=(a\otimes 1)\Delta(b)$ belong to $A\otimes A$ for any $a, b\in A$.
 \item[(2)] The counit satisfies $(\varepsilon\otimes id)T_{1}(a\otimes b) = ab = (id\otimes \varepsilon)T_{2}(a\otimes b)$.
 \item[(3)] $S$ is antimultiplicative and anticomultiplicative such that for any $a, b\in A$
 \begin{eqnarray}
 && S(a_{(1)}) a_{(2)(1)} \otimes a_{(2)(2)} = 1_{M(A)}\otimes a = a_{(1)} S(a_{(2)(1)}) \otimes a_{(2)(2)},  \label{2.1}\\
 && a_{(1)(1)} \otimes  S(a_{(1)(2)}) a_{(2)} = a\otimes 1_{M(A)} = a_{(1)(1)} \otimes a_{(1)(2)} S(a_{(2)}). \label{2.2}
 \end{eqnarray}
 \end{enumerate}
 If the antipode $S$ is bijective, then  multipler Hopf coquasigroup $(A, \Delta)$ is called \emph{regular}.
 Following the condition for antipode, we have
 \begin{eqnarray*}
 && m(id\otimes {S})\big( (a\otimes 1){\Delta}(b) \big) = {\varepsilon}(b)a, \\
 && m({S}\otimes id)\big( {\Delta}(a)(1 \otimes b) \big) = {\varepsilon}(a)b.
 \end{eqnarray*}

 A left (resp. right) integral on $A$ is a nonzero element $\varphi \in A^* := Hom(A, k)$ (resp. $\psi \in A^*$ ) such that
 \begin{eqnarray*}
 && (id\otimes \varphi)\Delta(a) = \varphi(a)1_{M(A)}  \quad \big(\mbox{resp.} (\psi\otimes id)\Delta(a) = \psi(a)1_{M(A)} \big), \quad \forall a\in A.
 \end{eqnarray*}
 Because $S$ is anticomultiplicative, we have that $\varphi\circ S$ is a right integral on $A$.

\section{Integrals on a multiplier Hopf coquasigroup}
\def\theequation{\thesection.\arabic{equation}}
\setcounter{equation}{0}

 Let $(A, \Delta)$ be a regular multiplier Hopf coquasigroup with a faithful integral $\varphi$.
 Just as in the case of algebraic quantum group (see Proposition 2.6 in \cite{DVZ}) or algebraic quantum hypergroups (see Proposition 1.6 in \cite{DV}),
 we show that the multiplier Hopf coquasigroup $(A, \Delta)$ must have local units in the sense of the following proposition.
 \\
 
 \textbf{Proposition \thesection.1}
 Let $(A, \Delta)$ be a regular multiplier Hopf coquasigroup with a non-zero integral $\varphi$.
 Given finite numbers of elements $\{a_{1}, a_{2}, \cdots, a_{n}\}$, there exists an element $e\in A$ such that $a_{i} e = a_{i} = e a_{i}$ for all $i$.

 \emph{Proof} It is similar to the proof of Proposition 1.6 in \cite{DV}. Set the linear space
 \begin{eqnarray*}
  V = \{(aa_{1}, aa_{2}, \cdots, aa_{n}, a_{1}a, a_{2}a, \cdots, a_{n}a) \mid a\in A\} \subseteq A^{2n}.
 \end{eqnarray*}
 Consider a linear functional on $A^{2n}$ that is zero on $V$. This means that we have functionals $w_{i}$ and $\rho_{i}$ on $A$ for $i = 1, 2, \cdots, n$, such that
 \begin{eqnarray*}
 \sum_{i=1}^{n} w_{i}(aa_{i}) + \sum_{i=1}^{n} \rho_{i}(a_{i}a) = 0 \quad \mbox{for all } a\in A.
 \end{eqnarray*}
 Then for all $x, a\in A$ we have
 \begin{eqnarray*}
 && x \Big(\sum_{i=1}^{n} (w_{i}\otimes id)\big(\Delta(a)(a_{i}\otimes 1)\big) + \sum_{i=1}^{n} (\rho_{i}\otimes id)\big((a_{i}\otimes 1)\Delta(a) \big)\Big) \\
 &=& \sum_{i=1}^{n} (w_{i}\otimes id)\big((1\otimes x)\Delta(a)(a_{i}\otimes 1)\big) + \sum_{i=1}^{n} (\rho_{i}\otimes id)\big((a_{i}\otimes 1)(1\otimes x)\Delta(a) \big)\\
 &=& 0,
 \end{eqnarray*}
 since $(1\otimes x)\Delta(a) \in A\otimes A$.
 Because the product in $A$ is nondegenerate, we have for all $a\in A$ that
 \begin{eqnarray*}
  \sum_{i=1}^{n} (w_{i}\otimes id)\big(\Delta(a)(a_{i}\otimes 1)\big) + \sum_{i=1}^{n} (\rho_{i}\otimes id)\big((a_{i}\otimes 1)\Delta(a) \big) = 0.
 \end{eqnarray*}
 Now applying $\varphi$ on this expression, we get
 \begin{eqnarray*}
 \varphi(a) \Big(\sum_{i=1}^{n} w_{i}(a_{i}) + \sum_{i=1}^{n} \rho_{i}(a_{i})\Big) = 0 \quad \mbox{for all } a\in A.
 \end{eqnarray*}
 As integral $\varphi$ is non-zero, we have
 \begin{eqnarray*}
 \sum_{i=1}^{n} w_{i}(a_{i}) + \sum_{i=1}^{n} \rho_{i}(a_{i}) = 0.
 \end{eqnarray*}
 So any lieaner functional on $A^{2n}$ that is zero on the space $V$ is also zero on the vector $(a_{1}, a_{2}, \cdots, a_{n}, a_{1}, a_{2}, \cdots, a_{n})$.
 Therefore $(a_{1}, a_{2}, \cdots, a_{n}, a_{1}, a_{2}, \cdots, a_{n}) \in V$.
 This means that there exists an element $e\in A$ such that $a_{i} e = a_{i} = e a_{i}$ for all $i$.
 $\hfill \Box$
 \\

 Recall from \cite{DV} that a linear functional $f$ on $A$ is called faithful if, for $a\in A$,
 we must have $a=0$ when either $f(ab) = 0$ for all $b\in A$ or $f(ba) = 0$ for all $b\in A$.
 Then under the faithfulness we can get the following result.
 \\
 
 \textbf{Lemma \thesection.2}
 Let $(A, \Delta)$ be a multiplier Hopf coquasigroup. If $f$ is a faithful linear functional on $A$,
 then for any $a\in A$ there is an element $e\in A$ such that
 \begin{eqnarray*}
 a = (id\otimes f)\big( \Delta(a)(1\otimes e) \big).
 \end{eqnarray*}

 \emph{Proof} Take $a\in A$ and set $V = \{ (id\otimes f)\big( \Delta(a)(1\otimes b) \big) \mid b\in A\}$, we need to show $a\in V$.
 Suppose that $a\notin V$, then there is on $A$ a functional $w\in A^*$ such that $w(a) \neq 0$ while $w|_{V} = 0$, i.e.,
 \begin{eqnarray*}
 0 = w\Big( (id\otimes f)\big( \Delta(a)(1\otimes b) \big) \Big) = f \Big( \big((w\otimes id) \Delta(a)\big) b \Big) \quad \mbox{for all } b\in A.
 \end{eqnarray*}

 Observe that $(w\otimes id) \Delta(a) \in M(A)$ and not necessarily belongs to $A$.
 However, we get $f \Big( \big((w\otimes id) \Delta(a)\big) b'b''\Big) = 0$ for all $b', b''\in A$ and by the faithfulness of $f$,
 we must have $\big((w\otimes id) \Delta(a)\big) b' =0$ for all $b'\in A$.

 If we apply the counit, $w(a)\varepsilon(b') = 0$ for all $b'\in A$ and hence $w(a) = 0$. This is a contradiction.
 $\hfill \Box$
 \\

 Similarly we have for a regular multiplier Hopf coquasigroup $(A, \Delta)$,
 \begin{eqnarray*}
 && a \in \{ (id\otimes f)\big( (1\otimes b)\Delta(a) \big) \mid b\in A\} \\
 && a \in \{ (f\otimes id)\big( (b\otimes 1)\Delta(a) \big) \mid b\in A\} \\
 && a \in \{ (f\otimes id)\big( \Delta(a)(b\otimes 1) \big) \mid b\in A\}
 \end{eqnarray*}
 for any faithful $f\in A^*$.
 In particularly, when we assume that a left integral $\varphi$ is faithful, then
 \begin{eqnarray*}
 && A = \mbox{span} \{ (id\otimes \varphi)\big(\Delta(a)(1\otimes b) \big) \mid a, b\in A\}, \\
 && A = \mbox{span} \{ (id\otimes \varphi)\big((1\otimes a)\Delta(b) \big) \mid a, b\in A\},
 \end{eqnarray*}
 where 'span' means the linear span of a set of element.
 \\

 Next, we give some equations on the left and right integral.

 \textbf{Proposition \thesection.3}
 Let $\varphi$ (resp. $\psi$) be a left (resp. right) integral on $A$, then for $a, b\in A$
 \begin{eqnarray}
 && a_{(1)} \varphi\big( a_{(2)}S(b) \big) = \varphi\big( aS(b_{(1)}) \big) b_{(2)}, \qquad a_{(1)}\varphi (ba_{(2)}) = S (b_{(1)}) \varphi (b_{(2)} a). \label{3.1} \\
 && \psi \big( S(a)b_{(1)}\big)b_{(2)} = \psi \big( S(a_{(2)})b \big) a_{(1)}, \qquad \psi(a_{(1)}b) a_{(2)} = \psi(a b_{(1)}) S(b_{(2)}). \label{3.2}
 \end{eqnarray}

 \emph{Proof}
 We prove the first two equations on $\varphi$, and the others are similar.
 \begin{eqnarray*}
  a_{(1)} \varphi\big( a_{(2)}S(b) \big)
  &\stackrel{(\ref{2.2})}{=}& a_{(1)} \big( S(b_{(1)(2)}) b_{(2)} \big) \varphi\big( a_{(2)} S(b_{(1)(2)}) \big) \\
  &=& \big(a_{(1)} S(b_{(1)(2)}) \big) b_{(2)} \varphi\big( a_{(2)} S(b_{(1)(2)}) \big) \\
  &=& \big(a_{(1)} S(b_{(1)})_{(1)} \big) b_{(2)} \varphi\big( a_{(2)} S(b_{(1)})_{(2)} \big) \\
  &=& \big(a S(b_{(1)}) \big)_{(1)} b_{(2)} \varphi\big( (a S(b_{(1)}))_{(2)} \big) \\
  &=& \varphi\big( aS(b_{(1)}) \big) b_{(2)},
 \end{eqnarray*}
 and
 \begin{eqnarray*}
  a_{(1)}\varphi (ba_{(2)})
  &\stackrel{(\ref{2.1})}{=}& \big( S(b_{(1)}) b_{(2)(1)} \big) a_{(1)} \varphi\big( b_{(2)(2)} a_{(2)} \big) \\
  &=& S(b_{(1)}) \big(b_{(2)(1)} a_{(1)}\big) \varphi\big( b_{(2)(2)} a_{(2)} \big) \\
  &=& S(b_{(1)}) \big(b_{(2)} a \big)_{(1)} \varphi\big( (b_{(2)} a)_{(2)} \big) \\
  &=& S (b_{(1)}) \varphi (b_{(2)} a).
 \end{eqnarray*}
 This completes the proof.
 $\hfill \Box$
 \\

 \textbf{Remark}
 (1) Following Lemma \thesection.2, we can easily check that (\ref{3.1}) and (\ref{3.2}) make sense.

 (2) These formulas are useful in the following part. We take the first one for example.
 When the antipode of $(A, \Delta)$ is bijective, $a_{(1)} \varphi\big( a_{(2)}S(b) \big) = \varphi\big( aS(b_{(1)}) \big) b_{(2)}$ is equivalent to
 \begin{eqnarray*}
 && S(id\otimes \varphi)\big(\Delta(a)(1\otimes b) \big) = (id\otimes \varphi)\big((1\otimes a)\Delta(b) \big),
 \end{eqnarray*}
 which is used to define the antipode in algebraic quantum hypergroup (see Definition 1.9 in \cite{DV}).

 Set $x = (id\otimes \varphi)\big(\Delta(a)(1\otimes b) \big)$ and apply $\varepsilon$ on the above equation,
 we have $\varepsilon\big( S(x) \big) = \varphi(ab) = \varepsilon(x)$. By Lemma \thesection.2 we get $\varepsilon \circ S = \varepsilon$.
 \\

 In the following, we will show the uniqueness of left faithful integrals.

 \textbf{Theorem \thesection.4}
 Let $\varphi'$ be another left faithful integral on $(A, \Delta)$, then $\varphi'=\lambda\varphi$ for some scalar $\lambda\in k$, i.e.,
 the left faithful integral on $A$ is unique up to scalar.

 \emph{Proof}
 From Proposition \thesection.3, we have $a_{(1)}\varphi (ba_{(2)}) = S(b_{(1)}) \varphi (b_{(2)} a)$ for all $a, b\in A$.
 Apply $\varphi'$ to both expressions in this equation. Because $\varphi' \circ S$ is a right integral, the right hand side will give
 \begin{eqnarray*}
 && \varphi'S (b_{(1)}) \varphi (b_{(2)} a) = \varphi (\varphi' S (b_{(1)})b_{(2)} a) = \varphi (\varphi' S (b) 1_{M(A)}\cdot a) = \varphi' S (b) \varphi (a).
 \end{eqnarray*}
 For the left hand side,
 \begin{eqnarray*}
 && \varphi'(a_{(1)}\varphi (ba_{(2)})) = \varphi'(a_{(1)})\varphi (ba_{(2)}) = \varphi (b\varphi'(a_{(1)})a_{(2)}) = \varphi (b \delta_{a}),
 \end{eqnarray*}
 where $\delta_{a} = \varphi'(a_{(1)})a_{(2)}$.
 Therefore, $\varphi' S (b) \varphi (a) = \varphi (b \delta_{a})$ for all $a, b\in H$.

 We claim that there is an element $\delta\in M(A)$ such that $\delta_{a} = \varphi(a)\delta$ for all $a\in A$.
 Indeed, for any $a'\in A$
 \begin{eqnarray*}
 \varphi (b \varphi (a')\delta_{a})
 &=& \varphi (a') \varphi (b \delta_{a}) = \varphi (a') \varphi'S(b) \varphi (a)\\
 &=& \varphi (a) \varphi'S(b) \varphi (a') = \varphi (a) \varphi (b \delta_{a'})\\
 &=& \varphi (b \varphi (a)\delta_{a'}),
 \end{eqnarray*}
 then $\varphi (a')\delta_{a} = \varphi (a)\delta_{a'}$ for all $a, a'\in A$, since $\varphi$ is faithful.
 Choose an $a'\in A$ such that $\varphi (a') = 1$ and denote $\delta = \delta_{a'}$, then $\delta_{a} = \varphi(a)\delta$.

 If we apply $\varepsilon$, we get
 \begin{eqnarray*}
 \varphi(a)\varepsilon(\delta) &=& \varepsilon(\delta_{a}) =  \varepsilon(\varphi'(a_{(1)})a_{(2)}) \\
 &=&  \varphi'(a_{(1)})\varepsilon(a_{(2)}) =  \varphi'(a_{(1)}\varepsilon(a_{(2)}) ) \\
 &=&  \varphi'(a)
 \end{eqnarray*}
 for all $a\in A$ and with $\lambda = \varepsilon(\delta)$, we find the desired result.
 $\hfill \Box$
 \\

 \textbf{Remark}
  (1) Similarly, the right faithful integral on $A$ is unique up to scalar.
  However, as in the special infinite dimensional Hopf algebra case, the non-zero faithful integrals do not always exist in infinite dimensional case.

  (2) The uniqueness of the faithful integral also provides the uniqueness of the antipode as in \cite{DV}.
 \\

 \textbf{Proposition \thesection.5}
 There is a unique invertible element $\delta\in M(A)$ such that for all $a\in A$
 \begin{enumerate}
 \item[(1)] $(\varphi\otimes id)\Delta(a) = \varphi(a)\delta$ \mbox{ and } $(id\otimes \psi)\Delta(a)= \psi(a)\delta^{-1}$.
 \item[(2)] $\varphi S(a) = \varphi(a\delta)$.
  \end{enumerate}

 \emph{Proof}
 In the proof of Theorem \thesection.4, $\varphi(a)\delta = \delta_{a} = \varphi'(a_{(1)})a_{(2)}$ and $\varphi' S (b) \varphi (a) = \varphi (b \delta_{a})$,
 we take $\varphi' = \varphi$ and get a element $\delta\in M(A)$ such that $(\varphi\otimes id)\Delta(a) = \varphi(a)\delta$ and $\varphi S (a) = \varphi (a\delta)$.
 This gives the first part of (1) and (2).

 If we apply $\varepsilon$ on the first equation, we find $\varepsilon(\delta)=1$.
 Because $S$ flips the coproduct and if we let $\psi=\varphi\circ S$, we get
 \begin{eqnarray*}
 (id\otimes \psi)\Delta(a)
 &=& S^{-1}(S\otimes \psi)\Delta(a) = S^{-1}(S\otimes \varphi\circ S)\Delta(a) \\
 &=& S^{-1}(id\otimes \varphi) (S\otimes S)\Delta(a) = S^{-1}(id\otimes \varphi) \Delta^{cop}(S(a)) \\
 &=& S^{-1}(\varphi\otimes id) \Delta(S(a)) \stackrel{(1)}{=} S^{-1}(\varphi(S(a))\delta) \\
 &=& \psi(a)S^{-1}(\delta).
 \end{eqnarray*}
 It remains to show $S^{-1}(\delta) = \delta^{-1}$.

 If we apply $\varphi$ to the formula (\ref{3.2}) $\psi(a_{(1)}b) a_{(2)} = \psi(a b_{(1)}) S(b_{(2)})$,
 we get
 \begin{eqnarray*}
 \psi(b) \varphi(a)
 &=& \psi(a b_{(1)}) \varphi S(b_{(2)}) = \psi(a b_{(1)}\psi(b_{(2)})) =  \psi(a \psi(b)S^{-1}(\delta))\\
 &=& \psi(b) \psi(aS^{-1}(\delta))
 \end{eqnarray*}
 for all $a, b\in A$. Then $\varphi(a) = \psi(aS^{-1}(\delta))$ for all $a\in A$.
 Therefore, $\varphi(a) = \varphi S(aS^{-1}(\delta)) = \varphi\big(aS^{-1}(\delta) \delta\big)$ and so $S^{-1}(\delta) \delta = 1_{M(A)}$.
 On the other hand, $\psi(a) = \varphi S(a) = \varphi (a\delta) = \psi(a\delta S^{-1}(\delta))$ and so $\delta S^{-1}(\delta) = 1_{M(A)}$.
 Hence, $\delta$ is invertible and $S^{-1}(\delta) = \delta^{-1}$, equivalently $S(\delta) = \delta^{-1}$.
 $\hfill \Box$
 \\

  \textbf{Remark} (1) The square $S^{2}$ leaves the coproduct invariant, it follows that the composition $\varphi\circ S^{2}$ of the left faithful integral $\varphi$ with $S^{2}$ will again a left faithful integral.
  By the uniqueness of left faithful integrals, there is a number $\tau \in k$ such that $\varphi\circ S^{2} = \tau\varphi$.

  (2) If we apply (2) in Proposition \thesection.3 twice, we get
 \begin{eqnarray*}
 && \varphi\big( S^{2}(a) \big) = \varphi\big( S(a)\delta \big) = \varphi\big( S(\delta^{-1} a) \big)  = \varphi\big((\delta^{-1} a) \delta \big) = \varphi\big(\delta^{-1} a \delta \big).
 \end{eqnarray*}
 So $\varphi\big(\delta^{-1} a \delta \big) = \tau\varphi(a)$.

 (3) We call $\delta$ the modular element as in algebraic quantum group. 
 Here we cannot conclude that $\Delta(\delta) = \delta \otimes \delta$ due to lack of the coassociativity of $\Delta$.
 \\

 Finally, just as in the algebraic quantum and algebraic quantum hypergroup case, we will show the existence of the modular automorphism.

 \textbf{Proposition \thesection.6}
 (1) There is a unique automorphism $\sigma$ of $A$ such that $\varphi(ab) = \varphi\big(b\sigma(a) \big)$ for all $a, b\in A$.
 We also have $\varphi\big( \sigma(a) \big) = \varphi(a)$ for all $a\in A$.

 (2) Similarly, there is a unique automorphism $\sigma'$ of $A$ satisfying $\psi(ab) = \psi\big(b\sigma'(a) \big)$ for all $a, b\in A$.
 And also $\psi\big( \sigma'(a) \big) = \psi(a)$ for all $a\in A$.

 \emph{Proof}
 (1) For any $p, q, x\in A$,
  \begin{eqnarray*}
 (\psi\otimes \varphi)\big(xq_{(1)} \otimes pS(q_{(2)})\big)
 &=& \psi(xq_{(1)}) \varphi(pS(q_{(2)})) = \varphi\big(p \underline{\psi(xq_{(1)})S(q_{(2)})} \big) \\
 &\stackrel{(\ref{3.2})}{=}& \varphi\big(p \psi(x_{(1)}q) x_{(2)} \big) =\psi(x_{(1)}q) \varphi\big(p x_{(2)} \big) \\
 &=& \psi\big(\underline{x_{(1)}\varphi(p x_{(2)})} q \big)  \stackrel{(\ref{3.1})}{=} \psi\big(S(p_{(1)})\varphi(p_{(2)} x) q \big) \\
 &=&\varphi(p_{(2)} x) \psi\big(S(p_{(1)}) q \big) = \varphi\Big((\psi\big(S(p_{(1)}) q \big)p_{(2)}) x \Big) \\
 &=& \varphi\Big((\psi\circ S \otimes id)\big( (S^{-1}(q) \otimes 1)\Delta(p) \big) x \Big).
 \end{eqnarray*}
 On the other hand, we also have
  \begin{eqnarray*}
 (\psi\otimes \varphi)\big(xq_{(1)} \otimes pS(q_{(2)})\big)
 &=& \psi(xq_{(1)}) \varphi(pS(q_{(2)})) = \psi\Big(x (q_{(1)}\varphi\big(pS(q_{(2)})\big)) \Big) \\
 &=& \psi\Big(x (id\otimes \varphi\circ S)\big(\Delta(q) (1\otimes S^{-1}(p))\big) \Big).
 \end{eqnarray*}
 Now assume that $\psi = \varphi \circ S$. then we have $\psi\circ S = \tau\varphi$ and $\psi(y) = \varphi(y\delta)$ by Proposition \thesection.5 (2).
 Then the above calculation will give us
 \begin{eqnarray*}
  \varphi(ax)
  &=& \psi(xb) = \frac{1}{\tau} \varphi\circ S(xb) \\
  &=& \frac{1}{\tau} \varphi(xb\delta) = \varphi\big( x(\frac{1}{\tau}b\delta) \big) \\
  &=& \varphi\big( xb \big)
 \end{eqnarray*}
 for all $x\in A$,  where $a= (\psi\circ S \otimes id)\big( (S^{-1}(q) \otimes 1)\Delta(p)$,
 $b' = (id\otimes \varphi\circ S)\big(\Delta(q) (1\otimes S^{-1}(p))\big)$ and $b = b'\delta$.

 Because $\varphi$ is faithful, the element $b$ is uniquely determined by the element $a$.
 So we can define $\sigma(a) = b$.
 Moreover, by Lemma \thesection.2 and its remark, all element in $A$ are of the form $a$ above,
 the map $\sigma$ is defined on all of $A$.
 The map $\sigma$ is injective by the faithfulness of $\varphi$,
 it is also surjective because all element in $A$ are also of the form $b$ above.

 Take $a, b, c\in A$, then
 \begin{eqnarray*}
  \varphi\big(c\sigma(ab)\big)
  &=& \varphi\big((ab)c\big) = \varphi(abc) = \varphi\big(a(bc)\big) \\
  &=& \varphi\Big((bc)\sigma(a)\Big) = \varphi\Big(b \big(c\sigma(a) \big)\Big) = \varphi\Big( \big(c\sigma(a) \big)\sigma(b)\Big)  \\
  &=& \varphi\Big(c \big(\sigma(a)\sigma(b) \big)\Big).
 \end{eqnarray*}
 It follows from the faithfulness of $\varphi$ that $\sigma(ab) = \sigma(a)\sigma(b)$. So $\sigma: A\longrightarrow A$ is an algebraic homomorphism.
 Apply this result twice, we have
 \begin{eqnarray*}
  \varphi(ab)
  = \varphi\big(b\sigma(a)\big)
  = \varphi\big(\sigma(a) \sigma(b)\big) =\varphi\big(\sigma(ab)\big) \quad \mbox{for all } a, b \in A.
 \end{eqnarray*}
 By Proposition \thesection.1 $A$ has local units, then $A^2 = A$, so $\varphi$ is $\sigma$-invariant.

 (2) Using that $\psi = \varphi \circ S$ we can easily get the statement for $\psi$.
 \begin{eqnarray*}
  \psi(ab) &=& \varphi S(ab) = \varphi\big( S(b)S(a) \big) \\
  &=& \varphi\Big( \sigma^{-1}S(a)S(b) \Big) = \varphi S \Big( b S^{-1}\sigma^{-1}S(a) \Big) \\
  &=& \psi \Big( b S^{-1}\sigma^{-1}S(a) \Big).
 \end{eqnarray*}
 Therefore, $\sigma' = S^{-1} \sigma^{-1} S$.
 $\hfill \Box$
 \\

 \textbf{Remark}
 In the proof of Proposition \thesection.6, we have
 \begin{eqnarray*}
 \varphi\Big((\psi\circ S \otimes id)\big( (S^{-1}(q) \otimes 1)\Delta(p) \big) x \Big)
 = \psi\Big(x (id\otimes \varphi\circ S)\big(\Delta(q) (1\otimes S^{-1}(p))\big) \Big).
 \end{eqnarray*}
 According to Lemma \thesection.2, we have that if $a\in A$ then there is a $b\in A$ such that $\varphi(ax)=\psi(xb)$ for all $x\in A$.
 This result will be used in the next section.
 \\

 As in the algebraic quantum group and hypergroup cases, the automorphism $\sigma$ and $\sigma'$ are called the modular automorphisms of $A$ associated with $\varphi$ and $\psi$ respectively.
 There are some extra properties derived from the above proposition.

 \textbf{Proposition \thesection.7} With the notation of above, we have
 \begin{enumerate}
 \item[(1)] $\sigma' = S^{-1} \sigma^{-1} S$ and $\sigma'(a) = \delta \sigma(a) \delta^{-1}$.
 \item[(2)] $\sigma(\delta) = \frac{1}{\tau} \delta$ and $\sigma'(\delta) = \frac{1}{\tau} \delta$
 \item[(3)] The modular automorphisms  $\sigma$ and $\sigma'$ commute with each other.
 \item[(4)] The modular automorphisms  $\sigma$ and $\sigma'$ commute with $S^2$.
 \item[(5)] For all $a\in A$, $\Delta\big(\sigma(a)\big) = (S^{2} \otimes \sigma)\Delta(a)$ and $\Delta\big(\sigma'(a)\big) = (\sigma'\otimes S^{-2})\Delta(a)$.
 \end{enumerate}

 Compared with algebraic quantum groups, $\Delta$ on  multiplier Hopf coquasigroup is not necessarily coassociatve.
  This is a significant difference between the two objects. Other than that, the proof is similar.

\section{Duality of discrete multiplier Hopf coquasigroups}
\def\theequation{\thesection.\arabic{equation}}
\setcounter{equation}{0}

 In this section, we will construct the dual of (infinite dimensional) multiplier Hopf coquasigroup of discrete type.
  The construction bases on the faithful integrals introduced in the last section.
 Here, we also start with defining the following subspace of the dual space $A^* $.
 \\

 \textbf{Definition \thesection.1}
 Let $\varphi$ be a left faithful integral on a regular multiplier Hopf coquasigroup $(A, \Delta)$.
 We define $\widehat{A}$ as the space of linear functionals on $A$ of the form $\varphi(\cdot a)$ where $a\in A$, i.e.,
 \begin{eqnarray*}
 \widehat{A} = \{\varphi(\cdot a) \mid a\in A\}.
 \end{eqnarray*}

 Because of Proposition 3.6 and the following remark, we have
 \begin{eqnarray*}
 && \widehat{A} = \{\varphi(\cdot a) \mid a\in A\} = \{\psi( a\cdot) \mid a\in A\} =  \{\varphi( a\cdot) \mid a\in A\} = \{\psi(\cdot a) \mid a\in A\}.
 \end{eqnarray*}

 Recall from \cite{Y} a regular multiplier Hopf coquasigroup $(A, \Delta)$ with a faithful integral $\varphi$ is called of discrete type,
 if there is a non-zero element $\xi \in A$ so that $a\xi=\varepsilon(a)\xi$ for all $a\in A$.

 The element $\xi$ ia called a \emph{left cointegral}. Similarly a \emph{right cointegral} is a non-zero element $\eta\in A$ so that $\eta a = \varepsilon(a)\eta$.
 The antipode will turn a left cointegral into a right one and a right one into a left one.

 Following this definition, we conclude $\varphi(\xi)\neq 0$. (If not, $0 = \varepsilon(a) \varphi(\xi) = \varphi(a\xi)$ for all $a\in A$,
 then $\xi=0$ by the faithfulness of $\varphi$. this is contradiction.)

 We start by making a discrete multiplier Hopf coquasigroup $(A, \Delta)$  into an unital algebra by dualizing the coproduct.
 \\

 \textbf{Proposition \thesection.2}
 For $w,w'\in \widehat{A}$, we can define a linear functional $ww'$ on $A$ by the formula
 \begin{eqnarray}
 && (ww')(x) = (w\otimes w')\Delta(x), \quad \forall x\in A. \label{4.1}
 \end{eqnarray}
 Then $ww'\in \widehat{A}$. This product on $\widehat{A}$ is not necessarily associative, but has a unit.

 \emph{Proof}
 Let $w, w'\in \widehat{A}$ and assume that $w' = \varphi(\cdot a)$ with $a\in A$. we have
 \begin{eqnarray*}
 (ww')(x) &=& (w\otimes \varphi(\cdot a))\Delta(x) = (w\otimes \varphi)\big(\Delta(x)(1\otimes a) \big) \\
 &=& w\big( x_{(1)} \varphi(x_{(2)}a) \big) \stackrel{(\ref{3.1})}{=} w \Big(S^{-1}\big( a_{(1)} \varphi(x a_{(2)}) \big)\Big) \\
 &=& \varphi\Big(x \big((wS^{-1}\otimes id)\Delta(a)\big) \Big)
 \end{eqnarray*}
 We see that the product $ww'$ is well-defined as a linear functional on $A$ and it has the form $\varphi(\cdot b)$, where $b = (wS^{-1}\otimes id)\Delta(a)$.
 So $ww' \in \widehat{A}$. Therefore, we have defined a product in $\widehat{A}$.

 The associativity of this product in $\widehat{A}$ is a consequence of the coassociativity of $\Delta$ on $A$, and
 $A$ is not necessarily coassociative.

 To prove that $\widehat{A}$ has a unit, assume that there is a cointegral $\xi \in A$ so that $a\xi=\varepsilon(a)\xi$ for all $a\in A$.
 \begin{eqnarray*}
 \varphi\big(\cdot \frac{1}{\varphi(\xi)} \xi\big) (a)
 &=& \frac{1}{\varphi(\xi)} \varphi\big(a \xi\big) =  \frac{1}{\varphi(\xi)} \varphi\big( \varepsilon(a) \xi\big) \\
 &=& \varepsilon(a),
 \end{eqnarray*}
 so $\varepsilon = \varphi\big(\cdot \frac{1}{\varphi(\xi)} \xi\big) \in \widehat{A}$.
 $\hfill \Box$
 \\

 \textbf{Remark}
 (1) Under the assumption, the elements of $\widehat{A}$ can be expressed in four different forms.
 When we use these different forms in the definition of product in $\widehat{A}$, we get the following useful expressions:
 \begin{enumerate}
 \item[(1)] $w\varphi(\cdot a) = \varphi(\cdot b)$  \mbox{with}  $b = wS^{-1}(a_{(1)}) a_{(2)}$;
 (2) $w\varphi( a\cdot) = \varphi(c \cdot)$  \mbox{with} $c = wS(a_{(1)}) a_{(2)}$.
 \item[(3)] $\psi(\cdot a)w = \psi(\cdot d)$  \mbox{with}  $d = a_{(1)} wS(a_{(2)})$; \quad
 (4)  $\psi( a\cdot)w = \psi(e \cdot)$  \mbox{with}  $e = a_{(1)} wS^{-1}(a_{(2)})$.
 \end{enumerate}

 (2) The reason for being restricted to the discrete case is that there is no definition of multiplier algebra $M(A)$ for a non-associative algebra $A$.
 \\

 Let us now define the comultiplication $\widehat{\Delta}$ on the unital algebra $\widehat{A}$.
 Roughly speaking, the coproduct is dual to the multiplication in $H$ in the sense that
 \begin{eqnarray*}
 \langle\widehat{\Delta}(w), x\otimes y\rangle = \langle w, xy\rangle, \quad \forall x, y\in H.
 \end{eqnarray*}

 We will first show that the above functional is well-defined and again in $\widehat{H}\otimes \widehat{H}$.

 \textbf{Proposition \thesection.3}
 Let $w\in \widehat{A}$, then we have $\widehat{\Delta}(w) \in \widehat{A} \otimes \widehat{A}$ and $\widehat{\Delta}$ is coassociative.

 \emph{Proof} The unit $1_{\widehat{A}} = \varepsilon = \psi\Big(\frac{1}{\psi S(\xi)} S(\xi) \cdot \Big) \in \widehat{A}$, and let $w = \psi(b\cdot)$. Then
 \begin{eqnarray*}
 \langle\widehat{\Delta}(w), x\otimes y\rangle
 &=& \langle (\varepsilon\otimes 1_{\widehat{A}})\widehat{\Delta}(w), x\otimes y\rangle  = \langle \varepsilon\otimes w, x_{(1)}\otimes x_{(2)}y\rangle\\
 &=& \langle \frac{1}{\psi S(\xi)}\psi\Big( S(\xi) \cdot \Big) \otimes \psi(b\cdot), x_{(1)}\otimes x_{(2)}y \rangle \\
 &=& \frac{1}{\psi S(\xi)}\psi\Big( S(\xi) x_{(1)} \Big) \psi(bx_{(2)}y) = \frac{1}{\psi S(\xi)} \psi\Big( b \underline{\psi\big( S(\xi) x_{(1)} \big)x_{(2)}} y \Big) \\
 &\stackrel{(\ref{3.2})}{=}& \frac{1}{\psi S(\xi)} \psi\Big( b \psi\big( S(\xi_{(2)}) x \big)\xi_{(1)} y \Big) = \psi\big( \frac{1}{\psi S(\xi)}S(\xi_{(2)}) x \big) \psi\Big( b\xi_{(1)} y \Big) \\
 &=& \langle \psi\big( \frac{1}{\psi S(\xi)}S(\xi_{(2)}) \cdot \big) \otimes \psi\Big( b\xi_{(1)} \cdot \Big), x\otimes y\rangle
 \end{eqnarray*}
 Hence $\widehat{\Delta}(w) = \psi\big( \frac{1}{\psi S(\xi)}S(\xi_{(2)}) \cdot \big) \otimes \psi\Big( b\xi_{(1)} \cdot \Big) \in \widehat{A} \otimes \widehat{A}$.

 The coassociativity is a direct consequence of the product associativity in $A$.
 $\hfill \Box$
 \\

 \textbf{Proposition \thesection.4}
 $\widehat{\Delta}: \widehat{A} \longrightarrow \widehat{A}\otimes \widehat{A}$ is an algebra homomorphism.

 \emph{Proof}
 It is straightforward that $\widehat{\Delta}$ is an algebra homomorphism, since for all $x, y\in H$
 \begin{eqnarray*}
 \langle\widehat{\Delta}(w_{1}w_{2}), x\otimes y\rangle &=& \langle w_{1}w_{2}, xy \rangle  = \langle w_{1}\otimes w_{2}, \Delta(xy) \rangle \\
 &=& \langle w_{1}, x_{(1)}y_{(1)} \rangle\langle w_{2}, x_{(2)}y_{(2)} \rangle \\
 \langle\widehat{\Delta}(w_{1})\widehat{\Delta}(w_{2}), x\otimes y\rangle
 &=&  \langle\widehat{\Delta}(w_{1})\otimes \widehat{\Delta}(w_{2}), (x_{(1)}\otimes y_{(1)}) \otimes (x_{(2)}\otimes y_{(2)}) \rangle \\
 &=& \langle\widehat{\Delta}(w_{1}), x_{(1)}\otimes y_{(1)} \rangle \langle\otimes \widehat{\Delta}(w_{2}), x_{(2)}\otimes y_{(2)} \rangle \\
 &=& \langle w_{1}, x_{(1)}y_{(1)} \rangle\langle w_{2}, x_{(2)}y_{(2)} \rangle
 \end{eqnarray*}
 This completes the proof.
 $\hfill \Box$
 \\

 Let $w\in \widehat{A}$ and assume $w=\varphi(\cdot a)$ with $a\in A$. Define $\widehat{\varepsilon}(w)=\varphi(a)=w(1_{M(A)})$.
 Then $\widehat{\varepsilon}$ is a counit on $(\widehat{A}, \widehat{\Delta})$ as follows.

 \textbf{Proposition \thesection.5}
 $\widehat{\varepsilon}: \widehat{A}\longrightarrow k$ is an algebra homomorphism satisfying
 \begin{eqnarray}
 (id\otimes \widehat{\varepsilon})\widehat{\Delta}\big( w \big) =  w = (\widehat{\varepsilon}\otimes id) \widehat{\Delta} \big( w \big)
 \end{eqnarray}
 for all $w \in \widehat{A}$.

 \emph{Proof}
 Firstly, let $w_{1} = \varphi(a\cdot)$ and $w_{2} = \varphi(b\cdot)$, then
 $w_{1}w_{2} = \varphi(c\cdot)$ with $c=\varphi\big( aS(b_{(1)})\big) b_{(2)}$.
 Therefore, if $\psi = \varphi \circ S$ we have
 \begin{eqnarray*}
 \widehat{\varepsilon}(w_{1}w_{2}) &=& \varphi(c) = \varphi\big( aS(b_{(1)})\big) \varphi(b_{(2)}) \\
 &=& \varphi\big( aS(b_{(1)}\varphi(b_{(2)}))\big) = \varphi(a) \varphi(b) \\
 &=& \widehat{\varepsilon}(w_{1})\widehat{\varepsilon}(w_{2}).
 \end{eqnarray*}

 Secondly, let $w = \varphi(\cdot a)$, then we have
 \begin{eqnarray*}
 \langle\widehat{\Delta}(w), x\otimes y\rangle
 &=& \langle\widehat{\Delta}(w)(1_{\widehat{A}}\otimes\varepsilon), x\otimes y\rangle  = \langle w\otimes \varepsilon, xy_{(1)}\otimes y_{(2)}\rangle\\
 &=& \langle \varphi(\cdot a) \otimes \varphi\big(\cdot \frac{1}{\varphi(\xi)} \xi\big), xy_{(1)}\otimes y_{(2)} \rangle \\
 &=& \frac{1}{\varphi(\xi)} \varphi(xy_{(1)} a)  \varphi\big(y_{(2)} \xi\big) =  \frac{1}{\varphi(\xi)} \varphi\Big(x \underline{y_{(1)}\varphi\big(y_{(2)} \xi\big)} a\Big) \\
 &\stackrel{(\ref{3.1})}{=}& \varphi\Big(x S^{-1}(\xi_{(1)})\varphi\big(y \xi_{(2)}\big) a\Big) \\
 &=& \langle  \varphi\Big(\cdot \frac{1}{\varphi(\xi)} S^{-1}(\xi_{(1)}) a\Big)\otimes \varphi\big(\cdot \xi_{(2)}\big), x\otimes y\rangle.
 \end{eqnarray*}
 Hence $\widehat{\Delta}(w) =  \varphi\Big(\cdot \frac{1}{\varphi(\xi)} S^{-1}(\xi_{(1)}) a\Big)\otimes \varphi\big(\cdot \xi_{(2)}\big)$.
 Therefore,
 \begin{eqnarray*}
 (id\otimes\widehat{\varepsilon}) \widehat{\Delta}(w)
 &=& \varphi\Big(\cdot \frac{1}{\varphi(\xi)} S^{-1}(\xi_{(1)}) a\Big) \varphi\big(\xi_{(2)}\big) \\
 &=& \varphi\Big(\cdot \frac{1}{\varphi(\xi)} S^{-1}\big(\xi_{(1)}\varphi(\xi_{(2)})\big) a\Big) =  \varphi\Big(\cdot  a\Big)\\
 &=& w.
 \end{eqnarray*}

% \begin{eqnarray*}
% (\widehat{\varepsilon}\otimes id) \widehat{\Delta}(w)
% &=& \varphi\Big(\frac{1}{\varphi(\xi)} S^{-1}(\xi_{(1)}) a\Big)\otimes \varphi\big(\cdot \xi_{(2)}\big)
% = \frac{1}{\varphi(\xi)} \varphi\Big(S^{-1}(\xi_{(1)}) a\Big)\otimes \varphi\big(\cdot \xi_{(2)}\big) \\
% &=& \frac{1}{\varphi(\xi)} \varphi\big(\cdot \varphi\Big(S^{-1}(\xi_{(1)}) a\Big)\xi_{(2)}\big)
% = \frac{1}{\varphi(\xi)} \varphi\Big(\cdot \varphi\Big(S^{-1}(\xi)_{(2)} a\Big)S(S^{-1}(\xi)_{(1)})\Big)\\
% &=& \frac{1}{\varphi(\xi)} \varphi\Big(\cdot a_{(1)} \varphi \big(S(\xi) a_{(2)}\big)\Big)
% =\frac{\varphi S(\xi)}{\varphi(\xi)} \varphi\Big(\cdot a\Big).
% \end{eqnarray*}

 Finally, from Proposition \thesection.3 $\widehat{\Delta}(\psi(b\cdot)) = \psi\big( \frac{1}{\psi S(\xi)}S(\xi_{(2)}) \cdot \big) \otimes \psi\Big( b\xi_{(1)} \cdot \Big)$.
 Then $(\widehat{\varepsilon}\otimes id) \widehat{\Delta}(w) = \psi\big( \frac{1}{\psi S(\xi)}S(\xi_{(2)}) \big) \psi\Big( b\xi_{(1)} \cdot \Big)
 = \frac{1}{\psi S(\xi)} \psi\Big( b\xi_{(1)}\psi\big(S(\xi_{(2)}) \big) \cdot \Big) = \psi(b\cdot)$.
 This completes the proof.
 $\hfill \Box$
 \\

 Let $\widehat{S}: \widehat{A}\longrightarrow \widehat{A}$ be the dual to the antipode of $A$, i.e., $\widehat{S}(w)=w\circ S$.
 Then it is easy to see that $\widehat{S}(w) \in \widehat{A}$, and we have the following property.

 \textbf{Proposition \thesection.6}
  $\widehat{S}$ is antimultiplicative and coantimultiplicative such that
   \begin{eqnarray*}
 && m(id\otimes m)(\widehat{S}\otimes id \otimes id)(\widehat{\Delta} \otimes id) = \widehat{\varepsilon}\otimes id = m(id\otimes m)(id\otimes \widehat{S} \otimes id)(\widehat{\Delta} \otimes id),\\
 && m(m\otimes id)(id\otimes \widehat{S} \otimes id)(id \otimes \widehat{\Delta}) = id\otimes \widehat{\varepsilon} =  m(m\otimes id)(id\otimes id \otimes \widehat{S})(id \otimes \widehat{\Delta}).
 \end{eqnarray*}

 \emph{Proof}
 For $w_{1}, w_{2}\in \widehat{H}$ and any $x\in H$,
 \begin{eqnarray*}
 \langle\widehat{S}(w_{1}w_{2}), x\rangle &=& \langle w_{1}w_{2}, S(x)\rangle = \langle w_{1}, S(x_{(2)})\rangle\langle w_{2}, S(x_{(1)})\rangle \\
 &=& \langle \widehat{S}(w_{1}), x_{(2)}\rangle\langle \widehat{S}(w_{2}), x_{(1)}\rangle = \langle \widehat{S}(w_{2})\widehat{S}(w_{1}), x\rangle
 \end{eqnarray*}
 This implies $\widehat{S}$ is antimultiplicative.
 \begin{eqnarray*}
 \langle\widehat{\Delta}\widehat{S}(w), x \otimes y\rangle
 &=& \langle \widehat{S}(w), xy \rangle = \langle w, S(xy) \rangle = \langle w, S(y)S(x) \rangle \\
 &=& \langle \widehat{\Delta}(w), S(y)\otimes S(x) \rangle = \langle \widehat{\Delta}^{cop}(w), S(x)\otimes S(y) \rangle \\
 &=& \langle(\widehat{S}\otimes\widehat{S}) \widehat{\Delta}^{cop}(w), x \otimes y\rangle,
 \end{eqnarray*}
 We conclude $\widehat{S}$ is coantimultiplicative.

 Finally, we show $m(id\otimes m)(\widehat{S}\otimes id \otimes id)(\widehat{\Delta} \otimes id) = \widehat{\varepsilon}\otimes id$,
  the other three formulas is similar.
 \begin{eqnarray*}
 && \langle m(id\otimes m)(\widehat{S}\otimes id \otimes id)(\widehat{\Delta} \otimes id) \big( w\otimes w' \big), x \rangle \\
 &=& \langle (\widehat{S}\otimes id \otimes id)(\widehat{\Delta} \otimes id) \big( w\otimes w' \big), (id\otimes \Delta)\Delta(x) \rangle \\
 &=& \langle(\widehat{\Delta} \otimes id) \big( w\otimes w' \big), (S\otimes id \otimes id)(id\otimes \Delta)\Delta(x) \rangle \\
 &=& \langle w\otimes w', (m \otimes id)(S\otimes id \otimes id)(id\otimes \Delta)\Delta(x) \rangle \\
 &=& \langle  w\otimes w', 1_{M(A)} \otimes x \rangle \\
 &=& \widehat{\varepsilon}(w) w'(x).
 \end{eqnarray*}
 This completes the proof.
 $\hfill \Box$
 \\

 From now, we get the first main result of this section.

 \textbf{Theorem \thesection.7}
 Let $(A, \Delta)$ be a regular multiplier Hopf coquasigroup of discrete type with a left faithful integral $\varphi$.
 Then $(\widehat{A}, \widehat{\Delta})$ is a Hopf quasigroup introduced in \cite{KM}.
 \\

 Let $\psi$ be a right faithful integral on $A$. For $w=\psi(a\cdot)$ we set $\widehat{\varphi}(w) = \varepsilon(a)$. Then we have the following result.

 \textbf{Proposition \thesection.8}
 The functional $\widehat{\varphi}$ defined above is a left faithful integral on Hopf quasigroup $(\widehat{A}, \widehat{\Delta})$.

 \emph{Proof}
 It is clear that $\widehat{\varphi}$ is non-zero. Assume $w = \psi(b\cdot)$, then by Proposition \thesection.3
 \begin{eqnarray*}
 \widehat{\Delta}(w) = \psi\big( \frac{1}{\psi S(\xi)}S(\xi_{(2)}) \cdot \big) \otimes \psi\Big( b\xi_{(1)} \cdot \Big).
 \end{eqnarray*}
 Therefore, we have
 \begin{eqnarray*}
 (id\otimes \widehat{\varphi})\widehat{\Delta}(w)
 &=& \psi\big( \frac{1}{\psi S(\xi)}S(\xi_{(2)}) \cdot \big) \varepsilon\Big( b\xi_{(1)} \Big) \\
 &=& \psi\big( \frac{1}{\psi S(\xi)}S(\xi) \cdot \big) \varepsilon(b)
 = \widehat{\varphi}(w) 1_{\widehat{A}}.
 \end{eqnarray*}

 Next, we show that $\widehat{\varphi}$  is  faithful.
 If $w_{1}, w_{2}\in \widehat{A}$ and assume $w_{1} = \psi(a\cdot)$ with $a\in A$, we have
 $w_{1} w_{2} =\psi \big(a_{(1)} w_{2}S^{-1}(a_{(2)}) \cdot \big)$.
 Therefore, $\widehat{\varphi}(w_{1} w_{2}) = w_{2}S^{-1}(a)$.
 If this is 0 for all $a\in H$, then $w_{2} = 0$,
 while if this is 0 for all $ w_{2}$ then $a=0$.
 This proves the faithfulness of $\widehat{\varphi}$.
 $\hfill \Box$
 \\

 If we set $\widehat{\psi} = \widehat{\varphi}\circ \widehat{S}$ as we do for the multiplier Hopf coquasigroup $(A, \Delta)$,
 we find that when $w = \varphi(\cdot a)$
 \begin{eqnarray*}
 \widehat{\psi}(w) = \widehat{\varphi}\circ \widehat{S}(w) =  \widehat{\varphi} (w\circ S) = \widehat{\varphi} \big( \varphi S(S^{-1}(a) \cdot)\big) = \varepsilon \big( S^{-1}(a) \big) = \varepsilon(a).
 \end{eqnarray*}

\section{Biduality}
\def\theequation{\thesection.\arabic{equation}}
\setcounter{equation}{0}

 Recall from \cite{Y} that the integral dual $\widehat{H}$ of an infinite dimensional Hopf quasigroup $H$ is a regular multiplier Hopf coquasigroup of discrete type.
 Specifically, let $H$ be an infinite dimensional Hopf quasigroup with a faithful left integral $\varphi$ and $\widehat{H} = \varphi(\cdot H)$, 
 if $\varphi(\cdot H) =  \varphi( H\cdot)$ and $\varphi\big((\cdot h)h'\big), \varphi\big(h'(h\cdot)\big) \in \widehat{H}$ for all  $h, h'\in H$,
 then $\widehat{H}$ is a regular multiplier Hopf coquasigroup of discrete type.
 
 And by Theorem 4.7 the integral dual $\widehat{\widehat{H}}$ of the regular multiplier Hopf coquasigroup of discrete type $\widehat{H}$ is Hopf quasigroup.
 Then, how is the relation between $H$ and $\widehat{\widehat{H}}$?
 Similarly, for a discrete multiplier Hopf coquasigroup $A$, $\widehat{A}$ is a Hopf quasigroup, the relation of $A$ and $\widehat{\widehat{A}}$ is what we care about.
 This is the content of the following theorem (biduality theorem).
 \\

 \textbf{Theorem \thesection.1} Let $(H, \Delta)$ be a Hopf quasigroup, and $(\widehat{H}, \widehat{\Delta})$ be the dual mutiplier Hopf coquasigroup of discrete type.
 For $h\in H$ and $f\in \widehat{H}$, we set $\Gamma(h)(f) = f(h)$. Then $\Gamma(h) \in \widehat{\widehat{H}}$ for all $h\in H$.
 Moreover, $\Gamma$ is an isomorphism between the Hopf quasigroups $(H, \Delta)$ and $(\widehat{\widehat{H}}, \widehat{\widehat{\Delta}})$.

 \emph{Proof} For $h\in H$, first we show that $\Gamma(h)$, as a linear functional on $\widehat{H}$, is in $\widehat{\widehat{H}}$. 
 Indeed, let $f = \varphi\big(\cdot S(h)\big)$ and take any $f'\in \widehat{H}$. 
 By Proposition 4.4 in \cite{Y}, $f'f = \varphi(\cdot h')$  where $h' = f'(h_{(2)}) S(h_{(1)})$.
 Therefore,  
 \begin{eqnarray*}
 && \widehat{\psi}(f'f) = \varepsilon(h') = f'(h) = \Gamma(h)(f').
 \end{eqnarray*}
 So $\Gamma(h) = \widehat{\psi}(\cdot f)$ and $\Gamma(h) \in \widehat{\widehat{H}} $.
 
 It is clear that $\Gamma$ is bijective between the linear space $H$ and $\widehat{\widehat{H}}$ because of the bijection of the antipode.
 $\Gamma$ respects the multiplication and comultiplication is straightforward because in both case the product is dual to the coproduct and vice versa.
 For details,
 \begin{eqnarray*}
 \langle \Gamma(hh'), f\rangle &=& \langle f, hh'\rangle = \langle \widehat{\Delta}(f), h\otimes h'\rangle \\
  &=& \langle (\Gamma\otimes\Gamma)(h\otimes h'), \widehat{\Delta}(f)\rangle = \langle\Gamma(h)\Gamma(h'), \widehat{\Delta}(f)\rangle,\\
 \langle \widehat{\widehat{\Delta}}\Gamma(h), f\otimes f'\rangle 
 &=& \langle \Gamma(h), ff'\rangle = \langle ff', h\rangle \\
 &=& \langle f\otimes f', \Delta(h)\rangle = \langle (\Gamma\otimes\Gamma)\Delta(h), f\otimes f'\rangle.
 \end{eqnarray*}
 Hence, $\Gamma$ is an isomorphism between $H$ and $\widehat{\widehat{H}}$.
 $\hfill \Box$
 \\
 
 Similarly, we can get another isomorphism. 

 \textbf{Theorem \thesection.2} Let $(A, \Delta)$ be a discrete multiplier Hopf coquasigroup, and $(\widehat{A}, \widehat{\Delta})$ be the dual Hopf quasigroup.
 For $a\in A$ and $w\in \widehat{A}$, we set $\Gamma(a)(w) = w(a)$. Then $\Gamma(a) \in \widehat{\widehat{A}}$ for all $a\in A$.
 Moreover, $\Gamma$ is an isomorphism between the multiplier Hopf coquasigroup $(A, \Delta)$ and $(\widehat{\widehat{A}}, \widehat{\widehat{\Delta}})$.
 \\
 
 As in the cases of algebraic quantum group and  algebraic quantum group hypergroup, all the results also hold for ( flexible (resp. alternative, Moufang)) multiplier Hopf ($*$-) coquasigroups.
 At the end of this section, we return to our motivating example of multipler Hopf coquasigroups.
 \\
 
 \textbf{Example \thesection.3}
 Let $G$ be a infinite (IP) quasigroup with identity element $e$, by definition $u^{-1}(uv) = v =(vu)u^{-1}$ for all $u, v\in G$.
 The quasigroup algebra $kG$ has a natrual Hopf quasigroup structure. $\delta_{e}$ is the left and right integral on $kG$. 
 The integral dual $k(G)$ introduced in \cite{Y} is a multipler Hopf coquasigroup of discrete type with the structure as follows.

 As an algebra, $k(G)$ is a nondegenerate algebra with the product
 \begin{eqnarray*}
 \delta_{u}\delta_{v} = \delta_{u, v}\delta_{v},
 \end{eqnarray*}
 and $1 = \sum_{u\in G} \delta_{u}$ is the unit in $M(k(G))$.
 The coproduct, counit and antipode are given by
 \begin{eqnarray*}
 \widehat{\Delta}(\delta_{u}) = \sum_{v\in G}\delta_{v}\otimes \delta_{v^{-1} u}, \quad \widehat{\varepsilon}(\delta_{u}) = \delta_{u, e},\quad \widehat{S}(\delta_{u}) = \delta_{u^{-1}}.
 \end{eqnarray*}
 The left integral $\widehat{\varphi}$ and right integral $\widehat{\psi}$ on $k(G)$ is the function that maps every $\delta_{u}$ to $1$.
 $\delta_{e}$ is the left and right cointegral in $k(G)$. 
 
 Now, we construct the dual of $k(G)$ as introduced in Section 4. Then 
 \begin{eqnarray*}
 \widehat{k(G)} = \{ \widehat{\varphi}(\cdot \delta_{u}) \mid u\in G\}.
 \end{eqnarray*}
 The element $\widehat{\varphi}(\cdot \delta_{u}) = \widehat{\psi}(\cdot \delta_{u})$ maps $\delta_{u}$ to $1$ and maps $\delta_{v} (v\neq u)$ to $0$.

 By Theorem \thesection.2, $kG \cong \widehat{k(G)}$ as Hopf quasigroups. The isomorphism $\Gamma: kG\longrightarrow \widehat{k(G)}$ is given by
 \begin{eqnarray*}
 \Gamma(u) &=& \widehat{\psi}\Big(\cdot \varphi\big(\cdot S(u) \big)\Big) =  \widehat{\psi}\Big(\cdot \delta_{e}\big(\cdot u^{-1} \big)\Big) \\
 &=& \widehat{\psi}\Big(\cdot \delta_{u}\Big).
 \end{eqnarray*}
 So if we identify $\widehat{\psi}(\cdot \delta_{u})$ with $u$, then $\widehat{k(G)} = kG$.
 
 By Theorem \thesection.3, $k(G) \cong \widehat{kG}$ as multiplier Hopf coquasigroups. The isomorphism $\Gamma: k(G)\longrightarrow \widehat{kG}$ is given by
 \begin{eqnarray*}
 \Gamma(\delta_{u}) &=& \widehat{\widehat{\psi}}\Big(\cdot \widehat{\varphi}\big(\cdot S(\delta_{u}) \big)\Big) =  \delta_{e}\Big(\cdot \widehat{\varphi}\big(\cdot \delta_{u^{-1}} \big)\Big) \\
 &=& \delta_{e}\Big(\cdot u^{-1}\Big) = u.
 \end{eqnarray*}
 So $\widehat{kG} = k(G)$.

\section*{Acknowledgements}

% The authors would like to thank the referee for his/her valuable comments.
 The work was partially supported by the China Postdoctoral Science Foundation (No. 2019M651764)
 and National Natural Science Foundation of China (No. 11601231).

% The first author, Tao Yang wishes to thank Prof. David Yetter, and the Kansas State University Department of Mathematics for hospitality during his stay in the United States as a visiting scholar.


\vskip 0.6cm
\begin{thebibliography}{30}

 \bibitem{DV} L. Delvaux and A. Van Daele (2011). Algebraic quantum hypergroups.
 \emph{Advance in Mathematics} 226: 1134-1167.

 \bibitem{DVZ} B. Drabant, A. Van Daele and Y. H. Zhang (1999). Actions of multiplier Hopf algebras.
 \emph{Communications in Algebra} 27(9): 4117-4172.

 \bibitem{K} J. Klim (2010). Integral theory for Hopf (co)quasigroups.
 \emph{arXiv}: 1004.3929.

 \bibitem{KM} J. Klim and S. Majid (2010). Hopf quasigroups and the algebraic 7-sphere.
 \emph{Journal of Algebra} 323: 3067-3110.

 \bibitem{V94} A. Van Daele (1994). Multiplier Hopf algebras.
 \emph{Transaction of the American Mathematical Society} 342(2): 917-932.

 \bibitem{V98} A. Van Daele (1998). An algebraic framework for group duality.
 \emph{Advance in Mathematics} 140(2): 323-366.

% \bibitem{V08} A. Van Daele (2008). Tools for working with multiplier Hopf algebras.
% \emph{Arabian Journal for Science and Engineering} 33(2C): 505-527.
%
% \bibitem{WW19} W. Wang and S. H. Wang (2019). Characterization of Hopf quasigroups.
% \emph{arXiv}: 1902.10141.

 \bibitem{Y} T. Yang (2020). Integral Dual of some infinite dimensional Hopf quasigroups.
 \emph{arXiv}: 2008.07199.



\end{thebibliography}
\end {document}